\documentstyle[12pt]{article}

\def\be{\begin{equation}}
\def\ee{\end{equation}}

\catcode`\@=11
\ifcase\@ptsize
 \font\tenmsa=msam10
 \font\sevenmsa=msam7
 \font\fivemsa=msam5
 \font\tenmsb=msbm10
 \font\sevenmsb=msbm7
 \font\fivemsb=msbm5
 \font\teneu=eufm10
 \font\seveneu=eufm7
 \font\fiveeu=eufm5
 \font\tenib=cmmib10
 \font\sevenib=cmmib7
 \font\fiveib=cmmib5
\or
 \font\tenmsa=msam10 scaled \magstephalf
 \font\sevenmsa=msam7 scaled \magstephalf
 \font\fivemsa=msam5 scaled \magstephalf
 \font\tenmsb=msbm10 scaled \magstephalf
 \font\sevenmsb=msbm7 scaled \magstephalf
 \font\fivemsb=msbm5  scaled \magstephalf
 \font\teneu=eufm10  scaled \magstephalf
 \font\seveneu=eufm7  scaled \magstephalf
 \font\fiveeu=eufm5   scaled \magstephalf
 \font\tenib=cmmib10  scaled \magstephalf
 \font\sevenib=cmmib7  scaled \magstephalf
 \font\fiveib=cmmib5   scaled \magstephalf
\or
 \font\tenmsa=msam10 scaled \magstep1
 \font\sevenmsa=msam7 scaled \magstep1
 \font\fivemsa=msam5  scaled \magstep1
 \font\tenmsb=msbm10 scaled \magstep1
 \font\sevenmsb=msbm7 scaled \magstep1
 \font\fivemsb=msbm5  scaled \magstep1
 \font\teneu=eufm10   scaled \magstep1
 \font\seveneu=eufm7 scaled \magstep1
 \font\fiveeu=eufm5 scaled \magstep1
 \font\tenib=cmmib10     scaled \magstep1
 \font\sevenib=cmmib7   scaled \magstep1
 \font\fiveib=cmmib5   scaled \magstep1
\fi
\newfam\msafam
\textfont\msafam=\tenmsa  \scriptfont\msafam=\sevenmsa
  \scriptscriptfont\msafam=\fivemsa
\newfam\msbfam
\textfont\msbfam=\tenmsb  \scriptfont\msbfam=\sevenmsb
  \scriptscriptfont\msbfam=\fivemsb
\def\Bbb{\ifmmode\let\next\Bbb@\else
 \def\next{\errmessage{Use \string\Bbb\space only in math mode}}\fi\next}
\def\Bbb@#1{{\Bbb@@{#1}}}
\def\Bbb@@#1{\fam\msbfam#1}
\newfam\eufam
\textfont\eufam=\teneu  \scriptfont\eufam=\seveneu
  \scriptscriptfont\eufam=\fiveeu
\def\frak{\ifmmode\let\next\frak@\else
 \def\next{\errmessage{Use \string\frak\space only in math mode}}\fi\next}
\def\frak@#1{{\frak@@{#1}}}
\def\frak@@#1{\fam\eufam#1}
\newfam\ibfam
\textfont\ibfam=\tenib  \scriptfont\ibfam=\sevenib
  \scriptscriptfont\ibfam=\fiveib
\def\bold{\ifmmode\let\next\bold@\else
 \def\next{\errmessage{Use \string\bold\space only in math mode}}\fi\next}
\def\bold@#1{{\bold@@{#1}}}
\def\bold@@#1{\fam\ibfam#1}

\begin{document}
\begin{titlepage}
\begin{flushright}
FTUV 96--38\\
IFIC 96--46\\
July,\, 1996
\end{flushright}
\vspace{2.5cm}
\begin{center} 
{\large\bf On the Completeness of Some Subsystems\\[2mm]
of $q$-deformed Coherent States}
\footnote{This is slightly modified version of the paper [Pe 1996]
 published in Helv. Phys. Acta 68, 554 (1996),}\\[2.5cm]
{\large\bf A. M. Perelomov}
\footnote
{On leave of absence from
Institute of Theoretical
and Experimental Physics, 117259 Moscow, Russia.
Email: perelomo@evalvx.ific.uv.es}\\
\medskip
{\em Departamento de F\'isica Te\'orica, Univ. de Valencia \\
46100-Burjassot (Valencia), Spain}\\[2.3cm]
{\bf Abstract}\\[0.5cm]
\end{center}
The von Neumann type subsystems of $q$-deformed coherent states are
considered. The completeness of such subsystems is proved.
\end{titlepage}

\section*{Introduction}
The systems of coherent states related to Lie groups introduced in [Pe 1972], 
play the important role in many branches of theoretical and mathematical 
physics and pure mathematics [CS 1985], [Pe 1986].

The basic feature of such systems is that they are overcomplete, i.e. 
contain subsystems, which are themselves complete. The most interesting of 
them are subsystems related to discrete subgroups of Lie groups, the first 
of which were considered by von Neumann [Ne 1929], [Ne 1932]. The completeness
properties of such system were investigated in [BBGK 1971] and 
[Pe 1971].

In the last few years $q$-deformed coherent states were introduced and some 
their properties were investigated (see, for example, [AC 1976], [Bi 1989], 
[Ma 1989], \mbox{[Ju 1991]).} Note that these states are related to 
$q$-deformed Lie algebras [KR 1981], [Dr 1985], [Dr 1986], [Ji 1985], 
[Ji 1986], [FRT 1991].

We investigate in the present paper the completeness properties of 
$q$-deformed coherent states for  simplest $q$-deformed Lie algebras, 
namely for $w_{q}(1), su_{q}(2)$ and $su_{q}(1, 1)$. It appears 
that some important properties of such systems are changed essentially after 
$q$-deformation.

\section{System of standard coherent states} 
\setcounter{equation}{0}
 
In this section we recall the basic properties of the system of standard
coherent states, introduced by E. Schr\"odinger [Sch 1926]
( see [Ste 1988] for historical discussion). 
For more details see books [CS 1985], [Pe 1986].
 
The basic quantities are the creation and annihilation operators $a^+$
 and $a$ and the unit operator $I$, which act in the Hilbert space ${\cal H}$
 and generate the Heisenberg--Weyl algebra: 
$$ [a, a^{+}] = aa^+ - a^{+}a = I,\quad [a, I] = [a^{+}, I] = 0.
\eqno (1.1)$$
 
 The standard orthonormal basis $\lbrace |n\rangle \rbrace ,\quad n = 0, 1,
 \ldots $, in  ${\cal H}$ is defined by
$$ |n \rangle = {(a^{+})^{n}\over \sqrt {n!}} |0\rangle ,
\eqno (1. 2)$$
 where $|0\rangle $ is the vacuum vector satisfying the  condition
 $$a\,|0\rangle = 0. \eqno (1. 3)$$
 
 The operators $a$ and $a^+$ act as follows
 $$a\,|n\rangle  = \sqrt {n}\,|n-1\rangle ,\quad a^{+}\,|\,n\rangle =
 \sqrt {n+1}\,|\,n+1\rangle .\eqno (1. 4)$$
 
 Let us introduce the operators
 $$E(\alpha ) = \exp\,(\overline {\alpha } a),\quad E^{+}(\alpha ) =
 \exp\,(\alpha a^{+}),\quad \alpha \in {\bf C}. \eqno (1. 5)$$
 
Then the standard system of coherent states that are non-normalized, may be 
defined by the formula
$$||\alpha \rangle = E^{+}(\alpha )\,|  0\rangle ,\eqno (1. 6)$$
or
$$||\alpha \rangle = \sum_{n=0}^{\infty }{\alpha ^{n}\over \sqrt {n!}}\,
|n\rangle . \eqno (1. 7)$$
 
It is easy to see that coherent states are eigenstates of the annihilation
operator
$$a\, || \alpha \rangle = \alpha \,|| \alpha \rangle ,\quad \alpha \in 
{\bf C},\eqno (1. 8)$$
and we can calculate the norm of such state
$$\langle \alpha || \alpha \rangle = \sum _{m,n=0}^{\infty } {\overline
{\alpha }^{m} \alpha ^{n}\over \sqrt {m!\,n!}}\, \langle m|n \rangle =
\sum _{n=0}^{\infty }\,{|\alpha |^{2n}\over \sqrt {n!}} = \exp\,(|\alpha
|^{2}). \eqno (1. 9)$$
Hence the normalized state $|\alpha \rangle $ has the form
$$|\alpha \rangle = \exp \,\biggl( -{|\alpha |^{2}\over 2}\biggr) ||\alpha 
\rangle =\exp \,\biggl( -{|\alpha |^2\over 2}\biggr) \sum {\alpha ^{n}\over 
\sqrt {n!}}\, |n\rangle. \eqno (1. 10)$$
The coherent states are not orthogonal to one another. The scalar product
 of two such states has the form
 $$\langle \alpha || \beta \rangle  = \exp \,(\bar {\alpha } \beta ).
 \eqno (1. 11)$$
 We also have the``resolution of the unity''
 $${1\over \pi }\int {d^{2}\,\alpha \,| \alpha \rangle \,\langle \alpha |} =
 \sum _{n=0}^{\infty }\,|n\rangle \,\langle n| = I,\eqno (1. 12)$$
 from which it follows that the system of coherent states is complete.
 
 This gives us the possibility to expand an arbitrary state $|\psi \rangle $ 
 on the states $|\alpha \rangle $
$$|\psi \rangle  = {1\over \pi } \int d^{2}\alpha \,c(\alpha )\,
|\alpha \rangle ,\quad c(\alpha )=\langle \alpha |\psi \rangle .\eqno(1. 13)$$
 
Note that if a coherent state $|\beta \rangle $ is taken as $|\psi
\rangle $, Eq. (1.13) defines a linear dependence between different
coherent states. It follows that the system of coherent states is
overcomplete, i.e. it contains subsystems that are complete.
 
Using (1.10) we obtain the following expression for $\langle \alpha |\psi
\rangle $ in (1.13):
$$\langle \alpha |\psi \rangle  = \exp \,\biggl( -{|\alpha |^{2}\over 2}
\biggr) \psi (\overline {\alpha }), \eqno (1. 14)$$
where
$$\psi (\alpha ) = \sum _{n=0}^{\infty }{c_n\over \sqrt {n!}}\, \alpha ^n,
\quad c_n = \langle n|\psi \rangle . \eqno (1. 15)$$
 
At the same time, the inequality $|c_{n}| = |\langle n|\psi \rangle |
\leq 1$ means that $\psi (\alpha )$ is an entire function of the complex 
variable $\alpha $ for the normalizing state $|\psi >$. 
We also have $|\langle \alpha |\psi \rangle | \leq 1$ and therefore 
have a bound on the growth of $\psi (\alpha )$:
$$|\psi (\alpha )|\leq \exp \,\biggl( {|\alpha |^{2}\over 2}\biggr) . 
\eqno (1. 16)$$
 
The normalization condition may now be written as
$$I = {1\over \pi } \int d^{2}\alpha \exp \,\biggl(-|\alpha |^{2}\biggr)\,
|\psi(\alpha )|^{2} = \langle \psi |\psi \rangle =1. \eqno (1. 17)$$
 
The expansion of an arbitrary state $|\psi \rangle $ with respect to
coherent states takes the form
$$|\psi \rangle  = {1\over \pi } \int d^{2}\alpha \exp \,\biggl( -{|\alpha |
^{2}\over 2}\biggr) \psi (\overline {\alpha })\,|\alpha \rangle . 
\eqno (1. 18)$$
 
Thus, we have established a one-to-one correspondence between the vectors
$|\psi \rangle $ of the Hilbert space and the entire functions
$\psi (\alpha )$, for which the integral (1.17) is finite. This
correspondence is established by Eqs. (1.15) and (1.18).

\section{System of coherent states for $q$-deformed Heisenberg-Weyl algebra} 

The generalization of the coherent states for $q$-deformed
Heisenberg--Weyl algebra was given in the papers [AC 1976], [Bi 1989],
[Ju 1991]. The corresponding formulae of the previous section 
should be modified.
 
Here the basic quantities as in the previous section are the creation and
annihilation operators $a^+$ and $a$ and unit operator $I$, which act in
the Hilbert space ${\cal H}$ and satisfy the relations
$$ [a, a^{+}] = aa^+ - qa^{+}a = I,\quad [a, I] = [a^{+}, I] = 0.
\eqno (2. 1)$$
 
 The orthonormal basis ${|n\rangle }$ in ${\cal H}$ is defined by 
 $$|n \rangle = {(a^{+})^{n}\over \sqrt {[n]!}}\, |0\rangle ,\eqno (2. 2)$$
 where
 $$[n]! = [1]\cdot [2]\cdot \ldots \cdot [n],\quad [n] = 1+q+\ldots
 +q^{n-1} = {1-q^n\over {1-q}}, \eqno (2. 3)$$
 and $|0\rangle $ is the vacuum vector satisfying the condition
 $$a\,|0\rangle = 0. \eqno (2. 4)$$
 
 The operators $a$ and $a^+$ act here as 
 $$a\,|n\rangle  = \sqrt {[n]}\,|n-1\rangle ,\quad a^{+}\,|n\rangle =
 \sqrt {[n+1]}\, | n+1\rangle . \eqno (2. 5)$$
 
 Let us introduce the operators
 $$E(\alpha ) = e_{q} (\overline {\alpha } a),\quad E^{+}(\alpha ) =
 e_{q} (\alpha a^{+}),\quad \alpha \in {\bf C},\eqno (2. 6)$$
where the function $e_{q}(x)$ is the generalization of the exponential
function and is defined by the formula (see [Ex 1983] and [An 1986] 
for details)\footnote{For simplicity, we restrict the 
consideration of the case $0 \leq q \leq 1$, used mainly in mathematical 
literature.}
 $$e_q(x) = \sum {x^n\over {[n]!}}. \eqno (2. 7)$$
It is easy to see that this series converges at $|x| < R_q =(1-q)^{-1}$
(for all finite values of $x$ at $|q|>1$), and
at $q \to 1,\quad [n]! \to n!$. This function coincides with a
standard $q$-exponent and satisfies the equation
$$\biggl( {d\over dx}\biggr) _{q} e_q (x) = e_q (x), \eqno (2. 8)$$
where the $q$-derivative $({d\over dx})_q$ is defined by the formula
$$\biggl( {d\over dx}\biggr) _{q}f(x) = {{f(x)-f(qx)}\over
{x(1-q)}}, \eqno (2. 9)$$
so that $({d\over dx})_{q} \to ({d\over dx})$ at the limit $q\to 1$. 
 
By using (2.8) and (2.9) one can show that
$$e_{q}(x) = {1\over \prod _{k=0}^{\infty }\,\biggl(1-q^k(1-q)x\biggr)}. 
\eqno (2. 10)$$
 
So the $e_q(x)$ is the meromorphic function, which has no zeros
and has simple poles at the points $x_k = q^{-k}/(1-q)$.
 
One can show [Ex 1983] that the inverse function $\biggl(e_{q}(x)
\biggr)^{-1}$ (an entire function) is given by 
$$\biggl(e_{q}(x)\biggr)^{-1} = e_{1/q}(-x) = \sum _{n=0}^{\infty }\,{{(-1)^
{n}q^{n(n-1)/2}\,x^{n}}\over {[n]!}} = \prod _{k=0}^{\infty }\,\biggl(1-q^
k(1-q)x\biggr). \eqno (2. 11)$$
 
We now define the system of coherent
states by the formula
$$||\alpha \rangle = E^{+}(\alpha )\, |0\rangle \eqno (2. 12)$$
and it is easy to see that this definition is coinside with the
definition given by Arik and Coon [AC 1976] 
$$||\alpha \rangle = \sum_{n=0}^{\infty }{\alpha ^{n}\over \sqrt {[n]!}}\,
|n\rangle . \eqno (2. 13)$$
 
It is easy to see that coherent states are eigenstates of the annihilation
operator
$$a\, || \alpha \rangle = \alpha \,|| \alpha \rangle ,\quad \alpha \in 
{\bf C},\eqno (2. 14)$$
and we can calculate the norm of such states
$$\langle \alpha || \alpha \rangle = \sum _{m,n=0}^{\infty } {\overline
{\alpha }^{m} \alpha ^{n}\over \sqrt {[m]!\,[n]!}} \,\langle m|n \rangle =
\sum _{n=0}^{\infty }{|\alpha |^{2n}\over  {[n]!}} = e_{q}\, \biggl(|\alpha
|^{2}\biggr). \eqno (2. 15)$$
Note that this series converges at
$$|\alpha |^{2} < R_{q}^2 = (1-q)^{-1}. \eqno (2. 16)$$
Hence the normalized state $|\alpha \rangle $ has the form
$$|\alpha \rangle = \biggl(e_{q} (|\alpha |^{2})\biggr)^{-1/2}\,||\alpha 
\rangle =\biggl(e_{q} (|\alpha |^2)\biggr)^{-{1\over 2}}\, \sum {\alpha ^{n}
\over \sqrt {[n]!}}\,|n\rangle ,\quad |\alpha | < R_{q}. \eqno (2. 17)$$
 
The coherent states are not orthogonal to one another. The scalar product
of two such states has the form
$$\langle \alpha || \beta \rangle  = e_{q}(\overline {\alpha }\beta ),
\eqno (2. 18)$$
We also have the ``resolution of the unity''
$${1\over \pi}\int _{D_{q}}d_q^{2}\alpha \,|\alpha \rangle \,\langle \alpha | 
={1\over 2\pi }\,\int _{0}^{2\pi }d\theta \int _{0}^{R_{q}^2} d_{q}(r^2)\,
|\alpha \rangle \,\langle \alpha | = \sum _{n=0}^{\infty }\,|n\rangle \,
\langle n|= I, \eqno (2. 19)$$
$$ \alpha =re^{i\theta },\quad D_q = \lbrace \alpha \colon |\alpha |<R_q 
\rbrace .$$
which follows from the formula
$$\int _{0}^{x_1}\, \biggl(e_{q}(x)\biggr)^{-1}x^{n}\, d_{q}x = [n]!, 
\eqno (2. 20)$$
where $x_{1} = (1-q)^{-1}$ is the first zero of the entire function
$\biggl(e_{q}(x)\biggr)^{-1}$ and the integral $\int _0^1 f(x)\, d_{q}x$ is 
the so-called Jackson integral [Ex 1983]:
$$\int _{0}^a f(x)\,d_{q}x = a(1-q)\, \sum _{k=0}^{\infty }\, q^{k}\, 
f(q^{k}a).\eqno (2. 21)$$
 
Note that from the``resolution of unity'' (2.19) it follows that the system
of coherent states is complete.
 This gives us the possibility to expand an arbitrary state $|\psi \rangle $ 
 on the states $|\alpha \rangle $
$$|\psi \rangle  = {1\over \pi }\, \int d^{2}\,\alpha \,\langle \alpha| 
\psi \rangle \,|\alpha \rangle . \eqno (2. 22)$$
 
If a coherent state $|\beta \rangle $ is taken as $|\psi \rangle $, (2.22) 
defines a linear dependence between different coherent states. Therefore, 
the system of coherent states is overcomplete, i.e.  contains subsystems, 
which are complete.
 
By using (2.17) we obtain an equation for $\langle \alpha |\psi \rangle $ 
$$\langle \alpha |\psi \rangle  = \biggl(e_{q} (|\alpha |^{2})\biggr)^{-1/2}
\,\psi(\overline {\alpha }), \eqno(2.23)$$
where
$$\psi (\alpha ) = \sum \,{c_n\over \sqrt {[n]!}}\, \alpha ^n,\quad c_{n} =
\langle n|\psi \rangle . \eqno(2.24)$$
 
At the same time, the inequality $|c_{n}| = |\langle n|\psi \rangle |
\leq 1$ means that $\psi (\alpha )$ for the normalizing state
$|\psi \rangle $ is an analytical function of the complex variable $\alpha $
in the disc $D_{q} = \{ \alpha || \alpha |  < R_{q} \} $.
We also have $|\langle \alpha |\psi \rangle | \leq 1$,
and therefore a bound on the growth of $\psi (\alpha )$:
$$|\psi (\alpha )| \leq \biggl(e_{q}(|\alpha |^{2})\biggr)^{1\over 2}. 
\eqno (2.25)$$
 
We can now rewrite the normalization condition as
$$I = {1\over \pi }\, \int _{D_q} d_q^{2}\alpha \,\biggl(e_{q}(|\alpha |^{2})
\biggr)^{-1}\,|\psi(\alpha )|^{2} = \langle \psi |\psi \rangle =1.
\eqno (2.26)$$
 
The expansion of an arbitrary state $|\psi \rangle $ with respect to
coherent states takes the form
$$|\psi \rangle  = {1\over \pi }\, \int _{D_q} d_q^{2}\alpha \,\biggl(e_{q}
(|\alpha |^{2})\biggr)^{-1/2}\,\psi (\overline {\alpha })\,|\alpha \rangle .
\eqno (2. 27)$$
 
Thus, we have established a one-to-one correspondence between the vectors
$|\psi \rangle $ of the Hilbert space and the functions $\psi (\alpha )$ 
analytical in $D_q$, for which the integral (2.26) is finite. This
correspondence is established by (2.23) and (2.27).

\section{Completeness of subsystems of $q$-deformed coherent states}
\setcounter{equation}{0}
 
As it was shown in the foregoing section, the system of $q$-deformed coherent 
states
$$\{|\alpha \rangle \colon \alpha \in D_{q}\},\quad D_{q} = \{\alpha
\colon |\alpha | \leq (1-q)^{-1/2}\}\eqno (3. 1)$$
is overcomplete, and hence there exist subsystems of coherent states 
that are complete ones. We describe these subsystems in this section. 
 
Let us take some set of points $\{\alpha _{k}\}$ in the disc $D_{q}$ and
take the corresponding subsystem of coherent states $\{|\alpha _{k}\rangle
\}$. Then if there exists a vector $|\psi \rangle $ of the Hilbert
space ${\cal H}$, which is orthogonal to all states $\{|\alpha _{k}
\rangle \}$:
$$\langle \alpha _{k}|\psi \rangle  = 0, \eqno (3. 2)$$
then the system $\{|\alpha _{k}\rangle \}$ is incomplete. It is complete
if such a vector does not exist.
 
We may reformulate this criterion in terms of the function
$$\psi (\alpha ) = \langle \psi || \alpha \rangle  = \sum _{n} \langle \psi
|n\rangle {\alpha ^{n}\over \sqrt {[n]!}}, \eqno (3. 3)$$
which is analytic inside $D_{q}$, and is equal to zero at the points
$\alpha _{k}$
$$\psi (\alpha _{k}) = 0. \eqno (3. 4)$$
If such a function has a finite norm
$$||\psi ||^{2} = \int _{D_q}|\psi (\alpha )|^{2}\, \biggl(e_{q}(|\alpha |^
{2})\biggr)^{-1}\, d_{q}^{2}\alpha  < \infty , \eqno (3. 5)$$
then the system $\{|\alpha _{k}\rangle \}$ is incomplete. But if any such
function has infinite norm $||\psi || = \infty $, then such a system is
complete.
 
Note that the function $\psi (\alpha )$, having the finite norm, should
satisfy the condition
$$|\psi (\alpha )|^{2}\, \biggl(e_{q}(|\alpha |^{2})\biggr)^{-1}\, \leq |\psi 
(\alpha )|^{2}\,\biggl(1-|\alpha |^{2}(1-q)\biggr) \leq C,\quad  \alpha 
\in D_{q}. \eqno (3. 6)$$
It follows from this condition that
$$\lim _{|\alpha |^{2} \to (1-q)^{-1}} |\psi (\alpha )|\, \biggl(1-|\alpha |^
{2}(1-q)\biggr)^{1\over 2} < \infty .\eqno (3. 7)$$
 
We give the simple example of the complete subsystem of coherent states.

Let the set $\{\alpha _{k}\}$ have a limit point inside the disc
$D_{q}$. The  function that is analytic inside $D_{q}$ and equal to zero at
points $\alpha _{k}$ should be equal to zero identically. Hence this system of
coherent states is complete.
 
For the future, it is convenient to introduce the new variable
$$\zeta  = (1-q)^{1/ 2}\,\alpha . \eqno (3. 8)$$
So we may now consider the set of functions analytical inside the unit disc
$D = \{\zeta \colon |\zeta | < 1\}$.
 
The characteristic property of $\psi (\zeta )$, related to the 
complete set $\{\zeta _{k}\}$, is that it has sufficiently many
zeros inside $D_{r} = \{\zeta \colon |\zeta | < r\}$ and hence it 
sufficiently quickly grows at $|\zeta | \to 1$. So we may use some theorems 
from the theory of functions analytical inside the unit disk. 

Let us give the theorem [Le 1964] that relates the growth of such function 
analytic in a disc with the distribution of its zeros.
 
Let $M(r)$ be the maximum modulus of  $f(\zeta )$ on the circle
$C_{r}=\lbrace \zeta \colon |\zeta | = r\rbrace $:
$$M_{1}(r) = \left[ {1\over 2\pi } \int |f(re^{i\theta })|^{2}\,
d\theta \right]^{1\over 2},
\eqno (3. 9)$$
and $n(r)$ be the number of zeros of $f(\zeta )$ in the disc $D=\lbrace 
\zeta \colon |\zeta |< r\rbrace $. We assume that the limit
$$\nu = \underline {\lim }_{r\to 1}\,(1-r^{2})\,n(r)\eqno (3. 10)$$
exists and that $\nu \neq 0$. The numbers
$$\tau  = \overline {\lim }_{r\to 1}\,\Big \lbrack \ln M(r)/\ln
{1\over {1-r^2}}\Big \rbrack ,\eqno (3. 11)$$
$$\tau _{1} = \overline {\lim }_{r \to 1}\,\Big \lbrack \ln M_{1}(r)/\ln
{1\over {1-r^2}}\Big \rbrack \eqno (3. 12)$$
characterize the growth of $f(\zeta )$ at $|\zeta |\to 1$, and 
we call $\tau $ and $\tau _{1}$ the generalized types of function $f(\zeta )$.
 
Note first of all the following relation between $\nu ,\,\tau $ and
$\tau _{1}$:

{\bf Theorem 3.1.}  {\em When $\nu  > 0$, the following inequalities are true}
$$\tau \geq {\nu \over 2},\quad \tau _{1} \geq {\nu \over 2}. \eqno (3. 13)$$

{\bf Proof.} Dividing $f(z)$ by $az^n$, if necessary, we obtain 
$\tilde {f}(z)$ with $\tilde {f}(0) = 1$. We use the Jensen formula
$${1\over 2\pi }\int _{0}^{2\pi }\ln |\tilde {f}(re^{i\theta })|\,d\theta  =
\int _{0}^{r} {n(t)\over t}\,dt, \eqno (3. 14)$$
from which 
$$\ln M(r) \geq \int _{0}^{r}{n(t)\over t}dt. \eqno (3. 15)$$
On the another hand, from the generalized inequality between the 
arithmetic and geometric means
$$\ln M_{1}(r) = {1\over 2}\ln\,\Big \lbrack {1\over 2\pi } \int _{0}^{2\pi }
|{\tilde f}(re^{i\theta })|^{2}\,d\theta \Big \rbrack  \geq {1\over 2\pi }
\int _{0}^{2\pi } \ln|{\tilde f}(re^{i\theta })|\,d\theta =\int _{0}^{r}{n(t)
\over t}\,dt.\eqno (3. 16)$$
 
It follows from the definitions of $\nu, \tau$ and $\tau _1$, that 
whatever the numbers $\varepsilon > 0,\quad \varepsilon _{1} > 0$
 and $\delta  > 0$, there exists the such number $r_{0} < 1$ that
 $$\ln M(r) \leq  (\tau +\varepsilon )\,\ln {1\over {1-r^2}},\quad \ln M_{1}
 (r) \leq (\tau _{1} + \varepsilon _{1})\,\ln {1\over {1-r^2}},\eqno (3. 17)$$
 $$n(r) \geq {{\nu  - \delta }\over {1-r^2}}\,r^2, $$
when $r>r_)$.
 We may now rewrite (3.15) and (3.16) as
 $$(\tau  + \varepsilon )\,\ln {1\over {1-r^2}} \geq \ln M(r) \geq \int _{0}^
 {r_0}{n(t)\over t}\,dt + {{\nu -\delta }\over 2}\,\Big \lbrack \ln
 {1\over {1-r^2}} - \ln {1\over {1-r_0^2}}\Big \rbrack, \eqno (3. 18)$$
 $$(\tau _{1} + \varepsilon _{1})\,\ln {1\over {1-r^2}} \geq \ln M_1(r) \geq
 \int _{0}^{r_0} {n(t)\over t}\,dt + {{\nu -\delta }\over 2}\,\Big \lbrack
 \ln {1\over {1-r^2}}- \ln {1\over {1-r_0^2}}\Big \rbrack. \eqno (3. 19)$$
 Considering the limit $r \to 1$ in (3.18) and (3.19) with $\nu > 0$, 
 we arrive at (3.13).
 
The criterion of completeness of the subsystem of $q$-deformed coherent 
states follows from this theorem and from inequality (3.6).

{\bf Theorem 3.2.}  {\it The system of $q$-deformed coherent states 
$\{|\alpha _{k}\rangle \}$ is complete, if the limit}
$$\nu  = \underline {\lim }_{r\to 1}\,(1-r^2)\,n(r)\eqno (3. 20)$$
{\it exists and if $\nu >1$. Here $n(r)$ is the number of points $\alpha _k$
inside the disc $D_{r}=\lbrace \zeta \colon |\zeta |< r\rbrace $.}
 
In order to construct the examples of complete subsystems it is useful to
consider the unit disc $D$ as a Lobachevsky plane with standard measure
$$d\mu (\zeta ) = {{d^{2}\zeta }\over {\bigl(1-|\zeta |^{2}\bigr)^{2}}},$$
on which the group $G=SU(1,1)/{\Bbb Z}_2$ acts transitively. The simplest
subsystems $\{|\alpha _{k}\rangle \}$ are related to the discrete subgroups
$\Gamma $ of the group $G$.
 
Let $\Gamma  = \{\gamma _{n}\}$ and $\alpha _0$ be any point of $D$.

{\bf Definition 3.3.} {\it The set of states $\{|\alpha _{k}
\rangle \}$, where $\alpha _{k} = \gamma _{k}\cdot \alpha _{0}$, is called 
the subsystem of coherent states related to subgroup $\Gamma $.}

 {\bf Theorem 3.4.} {\it The system of $q$-deformed coherent states related 
 to the discrete subgroup $\Gamma $ of the group $G = SU(1, 1)/{\Bbb Z}_{2}$ 
 is incomplete if the area $S_{\Gamma }$ of the fundamental domain $\Gamma
 \backslash D$ is infinite.}

 {\bf Proof.} In this case one may show (see for example [Le 1964]) that
 there exists a function $f(\zeta )$, which is analytic and bounded in $D$, 
 that has zeros at the points $\zeta _k$. For such a function the norm 
 defined by Eq. (3.5) is finite and hence this system of coherent states is
 incomplete.

 {\bf Theorem 3.5.} {\it Let the system of $q$-deformed coherent states 
 $\lbrace |\alpha _{k}\rangle \rbrace $ be related to the 
 discrete subgroup $\Gamma $, such that the area $S_{\Gamma }$ of the 
 fundamental domain $\Gamma \backslash D$ is finite and $S_{\Gamma }< \pi $.
 Then the system $\{|\alpha _{k}\rangle \}$ is complete.}

{\bf Proof.} Let us remind that non-Euclidian area of the disc of radius $r$ 
is equal to $S(r)=\pi r^2 /(1-r^{2})$. In this case, it follows from the 
condition $S_{\Gamma } < \pi $ that $\nu =\pi / S_{\Gamma } > 1$. 
Hence the norm of any analytic
function, which has the zeros in the points $\{\alpha _{k}\}$, is infinite,
and the system $\{|\alpha _{k}\rangle \}$ is complete.
 
Let us try to list the discrete subgroups $\Gamma $ for which $S_{\Gamma} <
\pi $. To this end, we need the information from the theory of discrete
 subgroups $\Gamma $ of the group $SU(1, 1)/{\Bbb Z}_{2}$, which we take from
 [Le 1964]. Let us restrict ourselves to consideration of groups
 with finite area of the fundamental domain $\Gamma \backslash D$. It is
 known that in this case the fundamental domain has the form of a polygon
 with an even number of sides $2n$. These sides being divided into pairs, 
 are equivalent with respect to the action of transformations of the group
 $\Gamma $. The  vertices of the polygon are joined in the cycles of
 vertices, which are equivalent to one another. With this, the sum of the
 angles of the polygon at the vertices of a given cycle equals to $2\pi /l$,
 where $l$ is either a positive integer or $\infty $. If $l=1$, the cycle
 is called random. If $l=\infty $, the vertices of the cycle lie on the
 boundary of the domain $D$, and the cycle is called parabolic, while in
 all the other cases, the cycle is called elliptic and $l$ is called
 the order of the cycle.  Let $c$ be the  number of cycles. By
 identifying equivalent sides and vertices, we obtain a Riemann surface. The
 genus  $p$ of this surface may be found by the formula
 $$2\pi = 1+n-c.\eqno (3. 21)$$
 
 We call the set of numbers $(p, c;\,l_1,\,l_2,\ldots ,\,l_c)$ the signature
 of the group $\Gamma $. We would like to mention that the area of the 
 fundamental domain
 $S_{\Gamma }$ is completely determined by the signature of the group and,
for our choice of invariant measure $d\mu (\zeta )=(1-|\zeta |^{2})^
 {-2}d\xi d\eta $, is given by 
 $$S_{\Gamma } = \pi \,\Big \lbrack p-1+{1\over 2}\sum _{j=1}^{c}\biggl( 1-{1
 \over l_{j}}\biggr) \Big \rbrack .\eqno (3. 22)$$
 
 From (3.22) it is easy to see that the value of $S_{\Gamma }$ cannot be
 arbitrarily close to zero. It may be shown [Si 1945] that
 the minimal value of $S_{\Gamma } = {\pi  \over 84}$ corresponds to the group
 $\Gamma $ with the signature  (0, 3;  2, 3, 7). If the fundamental domain
 is not compact, i.e. the group $\Gamma $ contains parabolic elements,
 then $S_{\Gamma } \geq {\pi \over 12}$; $S_{\Gamma } = {\pi \over 12}$
 correspond to the modular group $\Gamma = (0,3;\,2, 3, \infty )$. It is
  known also that when $p \geq 2$, the signature of $\Gamma $ may be
  arbitrary. For $p=1$ the condition $c \geq 1$ should be satisfied,
  and for $p=0$ we should have either $c \geq 5$, or $c=4$ and $\sum
  l_{j}^{-1} < 2$, or $c=3$ and $\sum l_{j}^{-1} < 1$.
 
  We are interested here in the case
  $${S_{\Gamma }\over \pi } = \Big\lbrack{p-1 + {1\over 2}\sum _{j=1}^{c}
  \biggl( 1-{1\over  l_j}\biggr) \Big\rbrack } < 1. \eqno (3. 23)$$

As it will shown below, the number of such cases is finite.

  Let us consider separately the different cases:
\begin{description} 
\item[I.]  Let $p \geq 2$, then  (3.23) cannot be satisfied.

\item[II.]  Let $p=1$, then  (3.23) takes the form
  $$\sum _{j=1}^{c} {1\over l_j} > c-2,\quad c \geq 1;\quad n=c+1.
\eqno (3.24)$$
Hence here we may have:
   \begin{description}
   \item[a)]
   $$c=1;\quad l_1 = 2, 3, \ldots , \infty ,\eqno (3.25)$$
   
   \item[b)]
   $$c=2;\quad l_1 =2, 3, \ldots , \infty ,\quad l_2 = 2, 3,\ldots , \infty ,
   \eqno(3.26)$$
    except of the case $l_{1} = \infty ,\, l_{2} = \infty $,
   
   \item[c)]
   $$c=3;\quad \sum _{j=1}^{3} {1\over l_{j}} > 1,\eqno (3.27)$$

   \begin{description}
   \item[$c_{1}$)]
   $$c=3;\quad l_{1} = 2,\,l_{2} =2,\,l_{3} =2, 3, \ldots ,< \infty ,
   \eqno (3.28)$$
   
   \item[$c_{2}$)]
   $$c=3;\quad (l_{1}, l_{2}, l_{3})=(2, 3, 3),\,(2, 3, 4),\,(2, 3, 5).
   \eqno (3.29)$$
   \end{description}
   \end{description}

\item[III.] Let
   \begin{description}
   \item[a)]
   $$p=0,\quad \sum _{j=1}^{c}{1\over l_{j}}>c-4,\eqno (3.30)$$

   \begin{description}
   \item[$a_{1}$)]
   $$c=5,\quad \sum _{1}^{5} {1\over l_{j}} > 1,\eqno (3.31)$$
   
   \item[$a_{2}$)]
   $$c=6,\quad \sum _{1}^{6} {1\over l_{j}} > 2,\eqno (3.32)$$
   
   \item[$a_{3}$)]
   $$c=7,\quad \sum _{1}^{7} {1\over l_{j}} > 3,\eqno (3.33)$$
   
   \item[$a_{4}$)]
   $$c=8,\quad \sum _{1}^{8} {1\over l_{j}} > 4.\eqno (3.34)$$
   This case and also the case $c > \infty $ are impossible.
   \end{description}
   
   \item[b)]
   $$c=4,\quad \sum _{1}^{4} {1\over l_{j}} > 0,\quad \sum _{j=1}^4
   {1\over l_{j}}<2,\eqno (3.35)$$
   
   \item[c)]
   $$c=3,\quad \sum _{1}^{3}{1\over l_{j}} > -1,\quad \sum _{j=1}^{3}
   {1\over l_{j}}<1.\eqno (3.36)$$
   \end{description}
\end{description}

So, as a function having zeros at the points $\zeta _n =\gamma _n\cdot 
\zeta _0$, 
we may take the automorphic form related to discrete subgroup $\Gamma =
\lbrace \gamma _n\rbrace $. If the fundamental domain $\Gamma \backslash D$ 
has 
finite area, we  may take it as polygon with finite number of sides which 
are segments of geodesics. Vertices of a polygon lying on the bondary of 
disc are called parabolic vertices. We denote ${\cal P}$ the set of parabolic 
vertices, and $D^+=D\bigcup {\cal P}$. Now we are ready to give 
the definition of the automorphic form.

{\bf Definition}. {\em An automorphic form of weight $m$ ($m$ is integer) is 
a function $f_{m}(z)$ that is analytic in $D$, satisfies the functional 
equation
\begin{displaymath} f_{m}(\zeta \cdot \gamma _{n})=(\beta _{n}z+\overline 
{\alpha }_{n})^{2m}\,f_{m}(z),\quad \gamma _{n} = \left( \begin{array}{cc}
\alpha _{n}&\beta _{n}\\
\overline{\beta }_{n}&\overline{\alpha }_{n}\end{array} \right) \in \Gamma 
\end{displaymath}
and is regular in $D^+$ (this means that at each parabolic vertex $\zeta _p$ 
of the domain $\Gamma \backslash D$, there should exist ${\rm lim}(z-z_{p})^
{2m}f_{m}(z)$ at $z\to z_p$, in the interior of domain $\Gamma \setminus D$). 
An automorphic form $f_{m}(z)$ is called parabolic if $f_{m}(z)$ vanishes at 
all parabolic vertices.}

The set of automorphic forms of weight $m$ builds a finite-dimensional vector 
space. We denote $d_{m}(\Gamma )$ ($d_{m}^{+}(\Gamma )$, correspondently) the 
dimension of the space of automorphic forms (the space of parabolic forms, 
correspodently). Let $m_{0}(m_{0}^+$) be the least $m$ for which $d_{m}
(\Gamma )\geq 2$ ($d_{m}^{+}(\Gamma )\geq 2$, correspondently). 
It is known (see for example [Le 1964]) that if $\Gamma \backslash D$ is 
compact, then $d_{m}(\Gamma )=d_{m}^{+}(\Gamma ),\,\,m_0 = m_{0}^{+}$, and 
any automorphic form may be considered as parabolic one.

The dimension of the space of automorphic forms of weight $m$ is given by
$$
d_{m}(\Gamma )=\cases{0,&for\,\,$m<0$,\cr
1,&for\,\,$m=0$,\cr
g_{1},&for\,\,$m=1$,\cr
(2m-1)(p-1)+\sum _{j=1}^{c}\,\left[ m\left( 1-\frac{1}{l_{j}}\right) \right],&
for\,\,$m\geq 2$.\cr}\eqno(3.37)$$

Here, $p$ is the genus of the fundamental domain, $[m]$ is the integer part 
of the number $m$, and $g_{1}\geq p$ is the number of holomorphic 
differentials on the Riemann surface $\Gamma \backslash D$.

With this, the number of zeros of function $f_{m}(z)$ in the interior of 
fundamental domain is given by the Poincar\'e formula [Po 1882] (written 
here in a somewhat different form)
$$N = 2mS_{\Gamma }/\pi .\eqno (3.38)$$

It should be mention that if there are elliptic and parabolic verticies, 
this number need not be integer.

Further, from a comparison of (3.22) and (3.37) we found that
$$N \geq  d_{m}+p-1,\eqno (3.39)$$
with the equality sign holding only in case when the numbers $m/l_{i}$ are 
integers, including zero.

In what follows we shall be interested in automorphic forms for which 
$d_{m}(\Gamma )\geq 2$. We denote by $m_{0}$ the minimal weight of such 
forms. We consider now the values which $m_{0}$ may take.
\begin{description}
\item[I.] If $p\geq 2$ then $m_{0}=1$, as it is evident from formula (3.37).
\item[II.] Let $p=1$ and let $c_{2}$ be the number of parabolic cycles. Then, 
if
\begin{description}
   \item[a)] $c_{2}\geq 2$, then $m_{0}=1$,
   \item[b)] $c_{2}=1$, then $m_{0}=2$,
   \item[c)] $c_{2}=0$ and $\Gamma =(1,1;\,2)$, then $m_{0}=4$,
   \item[d)] $c_{2}=0$ and $\Gamma =(1,1;\,l),\,l>3$, then $m_{0}=3$,
   \item[e)] $c_{2}=0, c\geq 2$, then $m_{0}=2$.
   \end{description}
\item[III.] If $p=0$ then
$$m_{0}\geq m_{1} = {\pi \over {2S_{\Gamma }}} = \biggl[\sum _{i=1}^{c}\,
\biggl(1-{1\over l_{i}}\biggr)-2\biggr]^{-1} = \biggl[c-2-\sum _{i=1}^{c}\,
{1\over l_{i}}\biggr]^{-1}.\eqno (3.40)$$
\end{description}

Let us introduce the notation:
$$N_{0}=2m_{0}\,S_{\Gamma }/\pi .$$ 
It follows from (3.39) that $N_{0}\geq p+1$.Therefore, $N_{0}$ can be equal 
to one only in the case of $p=0$.

With this, $d_{m_{0}}=2$, and the value of $m_{0}$ may be determined from 
(3.38):
$$m_{0} = {\pi \over 2S_{\Gamma }} = \biggl[\sum _{i=1}^{c}\,\biggl(1-{1\over 
l_{i}}\biggr)-2\biggr]^{-1}.\eqno (3.41)$$

Let $l$ be the least common multiple of the numbers $l_{j}$ which are not 
infinite. Then, (3.41) can be written in the form: 
$$m_{0}=l\biggl[\sum _{i}\,\biggl(l-{l\over l_{i}}\biggr)-2l\biggr]^{-1},$$ 
from which it follows that $m_{0}\leq 1$. However, $m_{0}$ must be divisible 
by those of $l_{i}$ which are not infinite. Therefore, they must coincide 
with $l$, i.e., $m_{0}=l$.

Thus we have

{\bf Proposition.} 
 If group $\Gamma $ of signature $(0,c;\,l_{1}\ldots 
l_{c_{1}},\,\infty ,\ldots \infty )$ admits automorphic form $f_{m_{0}}(z)$ 
with one zero in the fundamental domain, then $m_{0}$ is the least common 
multiple of the numbers $l_{1},l_{2},\ldots ,l_{c_{1}}$ and, moreover, must 
satisfy the condition
$$m_{0}(c-2)-\sum _{i=1}^{c_{1}}\,{m_{0}\over l_{i}} = 1.\eqno (3.42)$$

It is not difficult to show that (3.42) has a solution only for $c = 3,4,5$, 
and that the number of solutions of this equation is finite. There are 
21 discrete subgroups $\Gamma $ corresponding to them. All of them are 
listed in Table I.

\newpage
{\bf TABLE I}
\bigskip

\begin{tabular}{||r|l||}
\hline\hline
$m_{0}$\quad &\quad \quad \quad \quad $\Gamma $\\
\hline
1\quad &\quad (0,3;\,$\infty ,\infty ,\infty $)\\
2\quad &\quad (0,3;\,2,$\infty ,\infty $),\,(0,4;\,2,2,2,$\infty $),\,(0,5;\,
2,2,2,2,2)\\
3\quad &\quad (0,3;\,3,3,$\infty $)\\
4\quad &\quad (0,3;\,4,4,4),\,(0,3;\,2,4,$\infty $),\,(0,4;\,2,2,2,4)\\
6\quad &\quad (0,3;\,2,3,$\infty $),\,(0,3;\,3,3,6),\,(0,3;\,2,6,6),\,
(0,4;\,2,2,2,3)\\
8\quad &\quad (0,3;\,2,4,8)\\
10\quad &\quad (0,3;\,2,5,5)\\
12\quad &\quad (0,3;\,3,3,4),\,(0,3;\,2,3,12),\,(0,3;\,2,4,6)\\
18\quad &\quad (0,3;\,2,3,9)\\
20\quad &\quad (0,3;\,2,4,5)\\
24\quad &\quad (0,3;\,2,3,8)\\
42\quad &\quad (0,3;\,2,3,7)\\
\hline\hline
\end{tabular}
\bigskip

\section{Case of $q$-deformed $su(2)$-coherent states}
\setcounter{equation}{0}
 
In this section we consider the basic properties of the system of
$q$-deformed $su(2)$-coherent states (see [Ju 1991]
for other details).
 
The basic quantities here are the operators $J_{\pm }$ and $J_{0}$, which act
in the Hilbert space ${\cal H}$ of finite dimension $2j+1$ ($j$ is
half-integer, $2j+1$ is a positive integer) with the basis
$${|j, \mu \rangle },\quad \mu  =  -j, -j+1, \ldots ,j,\eqno(4.1)$$
or
$${|n\rangle },\quad n  =  j+\mu ,\quad n=0, 1, \ldots , 2j. \eqno(4.2)$$
The operators $J_{\pm }$ and $J_{0}$ act as follows
$$J_{\pm }\,|j, \mu \rangle  = {\sqrt {[j\mp \mu ]\,[j\pm \mu +1]}}\,
|j, \mu \pm 1\rangle ,\eqno(4.3)$$
$$J_{0}\,|j, \mu \rangle  = \mu \,|j, \mu \rangle ,\eqno(4.4)$$
or
$$J_{+}\,|n\rangle  =  \sqrt {[n+1]\,[2j-n]}\,|n+1\rangle ,\eqno(4.5)$$
$$J_{-}\,|n\rangle  =  \sqrt {[n]\,[2j-n+1]}\,|n-1\rangle ,\eqno(4.6)$$
$$J_{0} \,|n\rangle  =  (n-j) \,|n\rangle .\eqno(4.7)$$
Here $[n]$ is the Gauss symbol [Ga 1808]
$$ [n]=1+q+\ldots +q^{n-1} = {1-q^n\over 1-q},\quad [n]!=[1]\,[2]\ldots
[n].\eqno(4.8)$$
>From (4.5) it is not difficult to obtain
$$|n\rangle = \sqrt {{[2j-n]!\over [n]!\,[2j]!}}\,(J_{+})^{n}\,|0\rangle .
\eqno (4.9)$$
>From (4.5)--(4.7) it follows that the operators $J_{\pm }$ and $J_0$
satisfy the commutation relations
\footnote{Note that we use the quantum algebra generators different 
from the standard ones [KR 1981], [Ji 1985]}
$$[J_0, J_{\pm }] = \pm J_{\pm },\quad [J_{+}, J_{-}] = [2J_0],
\eqno (4. 10)$$
where the operator $[2J_0]$ is defined by the formula
$$[2J_0]\,|n\rangle = \lambda _{n}\,|n\rangle , \eqno (4. 11)$$
$$\lambda _{n} = \biggl(q^{j-\mu}[j+\mu ]-q^{j+\mu}[j-\mu ]\biggr) =
\cases{q^{j-\mu }\,[2\mu ],\quad \mu \geq 0;\cr
       -q^{j+\mu }\,[-2\mu ],\quad \mu \leq 0.\cr}\eqno (4.12)$$
 
Now we define the system of $q$-deformed 
coherent states by the formula
$$||z\rangle = e_{q}(zJ_{+})\,|0\rangle . \eqno (4. 13)$$
>From (4.9) we have
$$||z\rangle = \sum _{n=0}^{2j} \sqrt {{[2j]!\over [n]!\,[2j-n]!}}\,z^{n}
\,|n\rangle \eqno (4. 14)$$
and we may calculate the norm of this state
$$\langle z||z\rangle =G_{2j}\bigl(|z|^2\bigr)= \sum _{n=0}^{2j}\, {[2j]!\over
[n]!\,[2j-n]!}\,|z|^{2n} = \bigl[1+|z|^2\bigr]^{(2j)}. \eqno (4. 15)$$
Here $G_{2j}(x)$ is a certain polynomial of degree $2j$. Let us give the
simplest examples:
$$G_0 = 1,\quad G_1=1+x,\quad G_2=1+[2]x+x^2=1+(1+q)x+x^2;$$	
$$G_3 = 1+[3]x+[3]x^2+x^3 = (1+x)\,\biggl(1+([3]-1)x+x^2\biggr),\ldots 
\eqno (4. 16)$$
 
Note that these polynomials were first considered by Gauss [Ga 1808] and
investigated in more detail in the paper by Szeg\"o [Sz 1926]. Here we
note the following important properties of these polynomials:
\begin{description}
\item[i)]  Their roots are located on the circle of
unit radius and $x=1$ is not a root.

\item[ii)] The relation of these polynomials to theta--functions 
[Sz 1926]. Namely the functions
$$\Phi _{0}=1,\,\ldots ,\quad \Phi _{n} = {{(-1)^{n}\,q^{n/2}}\over {\sqrt
{(1-q)\,(1-q^{2})\ldots (1-q^{n})}}}\,G_{n}\biggl(-q^{-1/2}z\biggr)$$
are orthogonal on the unit circle $\lbrace z\colon \,z=e^{i\theta }\rbrace $
with the weight function $f(\theta )$, which coincides with theta--function
$$\int _{0}^{2\pi }\,\overline {\Phi }_{j}(\theta )\,\Phi _{k}(\theta )\,
f(\theta )\,d\theta  = 0,\quad j\neq k,$$
$$f(\theta ) = \sum _{n=-\infty }^{\infty }q^{n^{2}/2}\,e^{in\theta } = \sum
_{n=-\infty }^{\infty }q^{n^{2}/2}\cos n\theta =\biggl|D(e ^{i\theta })\biggr|
^{2},$$ 
$$\quad D(z) = \prod _{n=1}^{\infty }\,\sqrt {1-q^{n}}\,\biggl(1+q^{(2n-1)/2}z
\biggr).$$
\end{description}
Note also that we have
$$G_{n}\,\bigl(q^{1/2}\bigr)=\prod _{\nu =1}^{n}\,\bigl(1+q^{\nu /2}\bigr),$$
the expression for the generating function for $G_{n}(x)$,
$$\sum _{n=0}^{\infty }\,{{G_{n}(x)}\over {(1-q)\,(1-q^{2})\ldots (1-q^{n})}}
\,t^{n} = \prod _{n=0}^{\infty }{1\over {(1-q^{n}t)\,(1-q^{n}tx)}},$$
and the recurrence formulae
$$G_{n+1}(x) = (1+x)\,G_{n}(x)-(1-q^{n})\,xG_{n-1}(x);$$
$$G_{n}(qx)-(1-q^{n})\,G_{n-1}(qx) = q^{n}G_{n}(x).$$
 
Let us denote the roots as $\zeta _1,\ldots ,\zeta _{2j}.$ Then $|\zeta _k|=1$
and $\overline {\zeta }_k$ is also the root as $\zeta _k$. So
$$G_{n}(x)=\prod _{j=1}^{n}(x-\zeta _{j}).$$
The normalized coherent states now take the form
$$|z\rangle = \biggl(G_{2j}(|z|^2)\biggr)^{-1/2}\, \sum _{n=0}^{2j}\, 
\sqrt {{[2j]!\over [n]!\,[2j-n]!}}\,z^{n}\,|n\rangle ,\eqno (4. 17)$$
and the scalar product of two such states is
$$\langle w|z\rangle = {G_{2j}(\overline {w}z)\over \biggl(G_{2j}(|z|^2)\,
G_{2j}(|w|^2)\biggr)^{1/2}}. \eqno (4. 18)$$
So for the fixed coherent state $|z\rangle $ there are $2j$ coherent states
$|w_k\rangle ,\,k=1,\,\ldots ,\,2j$, which are orthogonal to state
$|z\rangle $. Here
$$w_k = (\overline {z})^{-1}\overline {\zeta }_k. \eqno (4. 19)$$
 
As for the standard system of coherent states, for $q$-coherent states we
also have the resolution of unity
$$\int ||z\rangle \,\langle z||\,d_{q}\mu (z) = I, \eqno (4. 20)$$
$$d_{q}\mu (z) = {[2j+1]\over 2\pi }\,\biggl(G_{2j+2}(|z|^2)\biggr)^{-1}\,
d_{q}(|z|^2)\,d\theta ,\quad z=|z|e^{i\theta }. \eqno (4. 21)$$
To prove this, let us consider the integral
$$I_{n, l}= \int _{0}^{\infty }\,x^{n}\,\biggl(G_{l}(x)\biggr)^{-1}\,d_{q}x. 
\eqno (4. 22)$$
Then after an integration by parts [Ex 1983] we have
$$I_{n,l} = {q^{-n}\,[n]\over [l-1]}\, \int _{0}^{\infty }\, x^{n-1}\,\biggl(
G_{l-1}(q^{-1}x)\biggr)^{-1}\,d_{q}x,\eqno (4.23)$$
and hence
$$I_{n,l} = {[n]!\over [l-1]\,[l-2]\ldots [l-n]}I_{0,l-n}. \eqno (4. 24)$$
Furthermore
$$I_{0,l-n}=\int _{0}^{\infty }\,\biggl(G_{l-n}(q^{-n}x)\biggr)^{-1}\,d_{q}x = 
{q^n\over [l-n-1]},\eqno (4. 25)$$
and finally
 
$$\int _{0}^{\infty }\,x^n\,\biggl(G_{l}(x)\biggr)^{-1}\,d_{q}x = {[n]!\,
[l-n-2]!\over [l-1]!}. \eqno (4. 26)$$
As a result of resolution of unity, an arbitrary vector $|\psi \rangle $
may be represented by a polynomial of degree $2j$:
$$\psi (\overline {z}) = \langle z|\psi \rangle . \eqno (4. 27)$$
Finally we come to the functional realization of the Hilbert space
${\cal F}_{j}$
$$\langle \psi _1|\psi _2\rangle = \int \overline {\psi _1(z)}\,\psi _2(z)\,
d_{q}\mu (z)\eqno (4. 28)$$
and have the basis:
$$f_{n}(\bar z) = \langle z||n \rangle ={ \sqrt {[2j]!\over
[n]!\,[2j-n]!}}\,\overline {z}^{n}. \eqno (4. 29)$$
It is easy to see that any set of (2j + 1) coherent states form
nonorthogonal basis in ${\cal F}_{j}.$

\section{Case of q-deformed $su_{q}(1,1)$-coherent states}
\setcounter{equation}{0}
 
In this section we consider the basic properties of the system of
$q$-coherent states for discrete series $T_k^+$  (see
[Ju 1991] for other details).
 
The basic quantities here are the operators $K_{\pm }$ and $K_{0}$, which act
in the infinite-dimensional Hilbert space ${\cal H}$ with the basis
$$\{|k,\, \mu \rangle \},\quad \mu =k,\, k+1,\, \ldots ,\eqno (5. 1)$$
or
$$\{|n\rangle \},\quad n=\mu -k,\, n=0,\,1,\, \ldots . \eqno (5. 2)$$
The operators $K_{\pm }$ and $K_{0}$ act as follows
$$K_{\pm }\,|k, \mu \rangle = {\sqrt {[\mu\pm k]\, [\mu \mp k\pm 1]}}\,|k,\,
\mu \pm 1\rangle , \eqno (5. 3)$$
$$K_{0}\,|k, \mu \rangle = \mu \,|k, \mu \rangle \eqno (5. 4)$$
or
$$K_{+}|n\rangle = \sqrt {[n+1]\,[2k+n]}\,|n+1\rangle , \eqno (5.5)$$
$$K_{-}|n\rangle =\sqrt {[n]\,[2k+n-1]}\,|n-1\rangle , \eqno (5.6)$$
$$K_{0}\, |n\rangle = (k+n)\, |n\rangle . \eqno (5.7)$$
Here $[n]$ is the Gauss symbol [Ga 1808]
$$[n]=1+q+\ldots +q^{n-1} = {1-q^n\over 1-q},\quad [n]!=[1]\,[2]\,\ldots
[n]. \eqno (5. 8)$$
>From (5.5) it is not difficult to obtain
$$|n\rangle = \sqrt {{[2k]!\over [n]!\,[2k+n-1]!}}\,(K_{+})^{n}\,|0\rangle .
\eqno (5.9)$$
>From (5.5)--(5.7) it follows that the operators $K_{\pm }$ and $K_0$
satisfy the commutation relations
\footnote{Note that we use the quantum algebra generators
different from the standard ones }
$$[K_0, K_{\pm }] = \pm K_{\pm },\quad [K_{-}, K_{+}] = [2K_0],
\eqno (5. 10)$$
where the operator $[2K_0]$ is defined by the formula
$$[2K_0]\,|n\rangle = \lambda _{n}\,|n\rangle , \eqno (5. 11)$$
$$\lambda _{n} = \biggl(q^{\mu -k}[\mu +k]+q^{\mu +k-1}[\mu -k ]\biggr) 
= \biggl(q^{2k+n-1}[n]+q^{n}[2k+n]\biggr) 
\eqno (5.12)$$
 
Now we define the system of $q$-coherent
states by the formula
$$||z\rangle = e_{q}\,(zK_{+})\,|0\rangle . \eqno (5. 13)$$
>From (5.9) we have
$$||z\rangle = \sum _{n=0}^{\infty} \,\sqrt {{[2k]!\over [n]!\,[2k+n]!}}\,z^n
\,|n\rangle \eqno (5. 14)$$
and may calculate the norm of this state
$$\langle z||z\rangle =F_{2k}\,\bigl(|z|^2\bigr)= \sum _{n=0}^{\infty}\, 
{[2k+n-1]!\over[n]!\,[2k]!}\,|z|^{2n} = \biggl(1-|z|^2\biggr)^{-(2k)}. 
\eqno (5. 15)$$
Here $F_{2k}(x)$ is the function of degree $(-2k)$:
$$F_{2k}(x)=G_{2k}^{-1}(-x).$$
Let us give the simplest examples:
$$F_0 = 1,\quad F_1=(1-x)^{-1},\quad F_2=\biggl(1-[2]x+x^2\biggr)^{-1}
=\biggl(1-(1+q)x+x^2\biggr)^{-1};$$	
$$F_3 =\biggl( 1-[3]x+[3]x^2-x^3\biggr)^{-1} = \Bigl((1-x)\,\biggl(1-([3]-1)x+
x^2\biggr)\Bigr)^{-1}\,\ldots \eqno (5. 16)$$
Here we only note that the poles of these functions are located on the
circle of unit radius and that, at integer $k$, $x=1$ is a pole.
 
Let us denote the poles as $\zeta _1,\ldots ,\zeta _{2k}.$ Then $|\zeta _k|=1$
and $\overline {\zeta }_k$ is the pole too as $\zeta _k$. So the normalized
coherent states have the form
$$|z\rangle = \biggl(F_{2k}(|z|^2)\biggr)^{-1/ 2}\, \sum _{n=0}^{\infty}\, 
\sqrt {{[2k]!\over [n]!\,[2k-n]!}}\,z^{n}\,|n\rangle ,\eqno (5. 17)$$
and the scalar product of two such states is
$$\langle w|z\rangle = {F_{2k}\,(\overline {w}z)\over \biggl(F_{2k}(|z|^2)\,
F_{2k}(|w|^2)\biggr)^{1/2}}. \eqno (5. 18)$$
 
As for the standard system of coherent states for $q$-coherent states we
also have the resolution of unity
$$\int ||z\rangle \,\langle z||\,d_{q}\mu (z) = I, \eqno (5. 19)$$
$$d_{q}\mu (z) = {[2k-1]\over 2\pi }\,\biggl(F_{2k+2}(|z|^2)\biggr)^{-1}\,
d_{q}(|z|^2)\,d\theta ,\quad z=|z|\,e^{i\theta }. \eqno (5. 20)$$
To prove this, let us consider the integral
$$I_{n, l}= \int _{0}^{1}\,x^{n}\,\biggl(F_{l}(x)\biggr)^{-1}\,d_{q}x. 
\eqno (5. 21)$$
Then after an integration by parts [Ex 1983] we have
$$I_{n,l} = {q^{-n}[n]\over [l-1]}\, \int _{0}^{1}\, x^{n-1}\,\biggl(F_{l-1}
(q^{-1}x)\biggr)^{-1}\,dx\eqno (5.22)$$
and hence
$$I_{n,l} = {[n]!\over [l-1]\,[l-2]\,\ldots [l-n]} I_{0,n-l}. \eqno (5. 23)$$
Furthermore
$$I_{0,n-l}=\int _{0}^{1}\biggl(F_{l-n}(q^{-n}x)\biggr)^{-1}\,d_{q}x 
= {q^n\over 
[l-n-1]}\eqno (5. 24)$$
and finally
$$\int _{0}^{1}\,x^n\,\biggl(F_{l}(x)\biggr)^{-1}\,d_{q}x = {[n]!\,[l-n-2]!
\over [l-1]!}.\eqno (5. 25)$$
As a result of resolution of unity, an arbitrary vector $|\psi \rangle $
may be represented by a function of degree $2k$:
$$\psi (\overline {z}) = \langle z||\psi \rangle . \eqno (5. 26)$$
And we finally  come to the functional realization of the Hilbert space
${\cal F}_{k}$:
$$\langle \psi _1|\psi _2\rangle = \int \overline {\psi _1(z)}\,\psi _2(z)\,
d_{q}\mu (z)\eqno (5. 27)$$
and we have the basis
$$f_{n}(\bar z) = \langle z||n \rangle ={ \sqrt {[2k]!\over
[n]!\,[2k-n]!}}\,\overline {z}^{n}.\eqno (5. 28)$$
So all formulae here are similar to the corresponding formulae for the case
of Heisenberg-Weyl algebra. Comparing, for example the basic formulae (2.15)
and (5.15) we can see that the case of Heisenberg-Weyl algebra is similar
to the case of $su$ (1,1) for $k={1\over 2}$.
So, for $su_{q}$(1,1) algebra we have the results
analogous to results of section 2.

\noindent
{\bf Acknowledgments.} It is pleasure to thank the Department of Theoretical
Physics, University of Valencia
for their hospitality and P. Kulish for useful remarks.

\end{document}